\documentclass[twoside]{article}
\title{Affine harmonic maps}
\author{J\"urgen Jost, Fatma Muazzez \c{S}im\c{s}ir}

\usepackage{a4,times,amsmath,theorem,latexsym,amssymb}
\usepackage{enumerate}

\newcommand{\be}{\begin{equation}}
\newcommand{\bel}{\begin{equation}}
\newcommand{\qe}{\end{equation}}
\newcommand{\ee}{\end{equation}}
\newcommand{\eeq}{\end{equation}}
\newcommand{\ba}{\begin{eqnarray}}
\newcommand{\ea}{\end{eqnarray}}
\newcommand{\rf}{\ref}

\newcommand{\bi}{\bibitem}

\newcommand{\R}{\mathbb{R}}

\newcommand{\Z}{\mathbb{Z}}

\newenvironment{summary}{\vskip\baselineskip \noindent\small\bf Summary: \rm}%
{\vskip\baselineskip}
\newenvironment{proof}{{\vskip\baselineskip\noindent\textbf{Proof:}}}%
{\hspace*{.1pt}\hspace*{\fill}\BOX\vskip\baselineskip}

\newcommand{\BOX}{\ensuremath\Box}
\newtheorem{theorem}{Theorem }[section]

\newtheorem{corollary}[theorem]{Corollary}

{\theorembodyfont{\rmfamily}}
{\theorembodyfont{\rmfamily}}
{\theorembodyfont{\rmfamily}}
\newtheorem{lemma}[theorem]{Lemma}

{\theorembodyfont{\rmfamily}}
{\theorembodyfont{\rmfamily}}

\begin{document}
\maketitle\thispagestyle{empty}


\begin{summary}
We introduce a class of maps from an affine flat into a Riemannian
manifold that solve an elliptic system defined by the natural second
order elliptic operator of the affine structure and the nonlinear
Riemannian geometry of the target. These maps are called affine
harmonic. We show an existence result for affine harmonic maps in a
given homotopy class when the target has nonpositive sectional
curvature and some global nontriviality condition is met. An example
shows that such a condition is necessary. \\
The analytical part is made difficult by the absence of a variational
structure underlying affine harmonic maps. We therefore need to
combine estimation techniques from geometric analysis and PDE theory
with global geometric considerations. 
\end{summary}

\section{Introduction}\label{section1}

A geometric structure usually induces a particular type of
connection that preserves that structure. When we have a Riemannian geometry, we get the Levi-Civita
connection as the unique torsion free connection that preserves the
metric. For a complex structure, we get a canonical complex connection. For
an affine structure -- which is the type of structure interesting us
in the present paper --, we obtain an affine flat connection. Thus, when we
have different  geometric structures on the same manifold, the induced
connections then in general are also different. For instance, for a
Hermitian metric on a complex manifold, its Levi-Civita connection
will in general not coincide with the holomorphic connection. More
precisely, the two coincide if and only if the manifold is
K\"ahler. This compatibility between two structures then makes the
theory of K\"ahler manifolds very rich. In fact, there is some analogy
between K\"ahler and  a particular class of affine structures first pointed out by Cheng and
Yau \cite{CY}. Remarkably, these structures also arise from a
completely different perspective, the one of information geometry,
that is, a geometric view of statistical families, see
\cite{Ch,AN,FA,J1}. \\

One of the motivations for the present work then is to develop
appropriate tools from geometric analysis to investigate such
structures. In Riemannian geometry, basic tools are geodesics and
harmonic maps. Here, for instance, a geodesic can be defined either
from a metric, as a curve that locally minimizes length, or from a
connection, as an autoparallel curve. The first one is a variational
characterization, the other is not. Likewise, harmonic maps are
characterized by a variational principle involving the metric. Since
harmonic maps are higher dimensional generalizations of geodesics, it
is then natural to develop also the corresponding concepts in terms of
a connection. This has been done by Jost-Yau \cite{JY} where the class
of Hermitian harmonic maps is introduced. These maps are determined by
the complex connection, and not by the Levi-Civita one. Therefore,
they do not satisfy a variational principle, and their investigation
becomes analytically much more difficult. Nevertheless, in \cite{JY},
a complete analysis could be carried out. As for ordinary harmonic
maps, it has to be required that the target manifold has nonpositive
sectional curvature. Still, an example in \cite{JY}  shows that in
contrast to ordinary harmonic maps, a Hermitian harmonic map need not
always exist in a given homotopy class, and a global nontriviality
condition needs to imposed to compensate for the lack of a variational
structure. \\

In this paper, we introduce the corresponding concept of affine
harmonic maps. They are determined in terms of an affine
connection. Thus, they also in general lack a variational
structure. In this paper, we succeed in extending the analysis of
\cite{JY} to affine harmonic maps and to show a general existence
theorem. \\

We hope that we can combine this existence results with Bochner type
identities in order to derive new restrictions on the topology of
affine flat manifolds.

\section{K\"ahler affine and dually flat manifolds}\label{section2}
We shall  use the standard conventions for raising and lowering indices.\\ 

An affine manifold $M$ possesses a covering by coordinate charts with affine coordinate changes. It then carries an affine flat connection, that is, one with vanishing curvature. This connection is complete if its geodesics can be defined on the real line. This is equivalent to the condition that the universal covering of $M$ is an affine vector space  which we identify with $\R^n$, with some abuse of notation. Note that compactness of $M$ does not imply its completeness. \\

It has been an important research topic to derive restriction on affine manifolds under various restrictions on their fundamental group, see e.g. \cite{Mil,Sm,GH,AyT}.\\

Cheng and Yau \cite{CY} then introduced an important condition which they called K\"ahler affine: $M$ carries a  2-tensor
\bel 
\label{1}
\gamma_{\alpha\beta} dx^\alpha dx^\beta
\qe
which locally is of the form
\bel
\label{2}
\gamma_{\alpha\beta}=\frac{\partial^2 F}{\partial x^\alpha \partial x^\beta}
\qe
for some convex function $F$, called a local potential (convexity here refers
to local coordinates $x$ and not to any metric.). Thus, $\gamma$ is positive definite and symmetric, that is, defines a Riemannian metric on $M$. In general, of course, the Levi-Civita connection of $\gamma$ will not be flat, that is, be different from the affine flat connection of $M$. The key point, however, is that the expression defining $\gamma$, 
\bel
\label{3}
\frac{\partial^2 F}{\partial x^\alpha \partial x^\beta} dx^\alpha dx^\beta
\qe
is invariant under affine transformations. \\
We can, however, recover the flat connection from $\gamma$ as follows: For $-1\le s \le 1$, we define
the $s $-connection through
\bel
\label{4}
\Gamma^{(s)}_{\alpha\beta\delta}= \Gamma^{(0)}_{\alpha\beta\delta}-\frac{s}{2}\partial_\alpha
\partial_\beta\partial_\delta F
\qe
where $\Gamma^{(0)}_{\alpha\beta\delta} $ represents  the Levi-Civita connection $\nabla^{(0)}$
 for $ \gamma_{\alpha\beta}$, i.e.,

\bel
\label{B50}
\Gamma^{(0)}_{\alpha\beta\delta}= \displaystyle \langle\nabla^{(0)}_{\frac{\partial}{\partial
x^\alpha}}\frac{\partial}{\partial x^\beta},
\frac{\partial}{\partial x^\delta}\rangle.
\qe
 Since by (\rf{2}),
\bel
\label{6}
\Gamma^{(0)}_{\alpha\beta\delta}=
\frac{1}{2}\partial_\alpha\partial_\beta\partial_\delta F,
\qe
we have
\bel
\label{7}
\Gamma^{(s)}_{\alpha\beta\delta}= \frac{1}{2}(1-s)\partial_\alpha
\partial_\beta\partial_\delta F,
\qe
and since this is
symmetric in $\alpha $ and $\beta $, $\nabla^{(s)} $ is torsion free.  Since
$\Gamma^{(s)}_{\alpha\beta\delta}+\Gamma^{(-s)}_{\alpha\beta\delta}=2\Gamma^{(0)}_{\alpha\beta\delta}$,
$ \nabla^{(s)}$ and
$\nabla^{(-s)} $ are dual to each other, in the sense that
\begin{equation}
Z\langle V,W\rangle=\langle\nabla^{(s)}_ZV,W\rangle+\langle V,\nabla^{(-s)}_ZW
\rangle
\end{equation}
for all vector fields $ V$, $W $, $Z $ where $\langle.,.\rangle$ stands for the metric $g$.\\ 
In particular, $
\Gamma^{(1)}_{\alpha\beta\delta}=0$, and so
$\nabla^{(1)} $ defines
a flat structure, and the coordinates $x$ are affine coordinates
for $\nabla^{(1)}$.\\
 The
connection dual to  $\nabla^{(1)}$ then is $\nabla^{(-1)} $ with Christoffel symbols
\[
\Gamma^{(-1)}_{\alpha\beta\delta}= \partial_\alpha
\partial_\beta\partial_\delta F\\
\]
with
respect to the $x$- coordinates. We can then obtain dually affine coordinates
$\xi$ by
\bel
\label{con7}
\xi_\beta=\partial_\beta F,\\
\end{equation}
and so also
\bel
\gamma_{\alpha\beta}=\partial_\alpha \xi_\beta.\\
\end{equation}
The
corresponding local potential is  obtained by a Legendre
transformation
\begin{equation}
\label{con9}
\Phi(\xi)=\max_{x}(x^\alpha\xi_\alpha-F(x)),\quad
F(x)+\Phi(\xi)-
x\cdot\xi=0,\\
\end{equation}
and
\begin{equation}
\label{con10}
x^\beta=\partial^\beta\Phi(\xi),\quad \gamma^{\alpha\beta}=\frac{\partial
x^\beta}{\partial \xi_\alpha}
=\partial^\alpha\partial^\beta\Phi(\xi).\\
\end{equation}

Thus, a K\"ahler affine structure yields a dually flat structure, i.e., a
Riemannian metric $ \gamma$ together with two flat connections $\nabla $ and
$\nabla^{*}$ that are dual with respect to $\gamma $. Such dually flat
structures have been introduced and investigated by Chensov \cite{Ch}
and Amari (see \cite{AN,FA}) as the basis of information
geometry. Conversely, given such a dually flat structures, one finds
local potential functions, that is, obtains a K\"ahler affine
structure, see e.g. the exposition in \cite{J1}. Thus, the two types
of structure are equivalent. Here, we work with the notion of K\"ahler
affine structure of Cheng-Yau because it is geometrically simpler and
more transparent.\\\\
Throughout this paper, we shall use standard summation
conventions. $c$ will denote a constant in estimates, without implying
that $c$ always has the same value. Being a constant here means that
it depends only on the underlying geometries as well as possibly on
the initial values, but not on the solutions
of the differential equations under consideration.

\section{Affine harmonic maps}\label{section3}
K\"ahler affine structure (\rf{2})  allows us to define a differential operator, 
\bel
\label{11} 
L:=\gamma^{\alpha\beta}\frac{\partial^2}{\partial x^\alpha\partial x^\beta},
\qe
that is affinely invariant. A function$f:M \to \R$ that satisfies
\bel
\label{12}
Lf=0
\qe
is called {\it affine harmonic}. More generally, when $N$ is a Riemannian manifold with metric $g_{ij}$ and Christoffel symbols $\Gamma^{i}_{jk}$, we call a map $f:M\to N$ {\it affine harmonic} if it satisfies
\bel
\label{13}
\gamma^{\alpha\beta}(\frac{\partial^2 f^i}{\partial x^\alpha\partial x^\beta}+\Gamma^{i}_{j k}\frac{\partial f^j}{\partial x^\alpha}\frac{\partial f^k}{\partial x^\beta})=0
\qe
in local coordinates on $N$. More invariantly, we can write (\rf{13}) as
\bel
\label{14}
\gamma^{\alpha\beta} D_\alpha D_\beta f=0
\qe
where $D$ is the connection on $T^\ast M \otimes f^{-1}TN$ induced by the flat connection on $M$ and the Levi-Civita connection on $N$. \\
We have the following general existence result for affine harmonic maps.
\begin{theorem}
Let $M$ be a compact K\"ahler affine manifold, $N$ a compact Riemannian manifold of nonpositive sectional curvature. Let $g:M\to N$ be continuous, and suppose $g$ is not homotopic to a map $g_0:M\to N$ for which there is a nontrivial parallel section of $g_0^{-1}TN$.\\
Then $g$ is homotopic to an affine harmonic map $f:M\to N$.
\end{theorem}
Using the argument first introduced by Al'ber \cite{Al2},
one can also show that the affine harmonic map is unique in its
homotopy class under the conditions of our  theorem.\\
After stating some corollaries and discussing an example, we shall
obtain this result in the next section by the method of \cite{JY}.\\

\begin{corollary}
Let $M$ be a compact K\"ahler affine manifold, $N$ a compact Riemannian manifold of negative sectional curvature. Let $g:M\to N$ be continuous, and suppose $g$ is not homotopic to a map onto a closed geodesic of $N$. Then $g$ is homotopic to an affine harmonic map.
\end{corollary}
\begin{corollary}
Let $M$ be a compact K\"ahler affine manifold, $N$ a compact Riemannian manifold of nonpositive sectional curvature. Let $g:M\to N$ be smooth and satisfy $e(g^\ast TN)\ne 0$, where $e$ is the Euler class. Then $g$ is  homotopic to an affine harmonic map.
\end{corollary}
The two corollaries follow from the theorem because their assumptions
imply that  $g$ cannot be homotopic to a map $g_0:M\to N$ for which
there is a nontrivial parallel section of $g_0^{-1}TN$. In fact, for
the first corollary, we observe that if the tangent space of $g_0(M)$
possesses a parallel section then $g_0(M)$ itself has to be a flat
subspace of the nonpositively curved space $N$. Since $N$ here is
assumed to have negative curvature, the only such subspaces are
one-dimensional, and they are homotopic to closed geodesics. For the
second corollary, we observe that a vector bundle with a parallel
section has vanishing Euler class. \\

(\rf{13}) is a semilinear system of elliptic partial differential equations. It is in general not in divergence form, and therefore, variational methods are not available for its investigation. 
The method of \cite{JY} which we shall use for these existence theorems consists in studying the associated parabolic equation,
\bel
\label{15}
\frac{\partial f^i}{\partial t}=\gamma^{\alpha\beta}\left(\frac{\partial^2 f^i}{\partial x^\alpha\partial x^\beta}+\Gamma^{i}_{jk}\frac{\partial f^j}{\partial x^\alpha}\frac{\partial f^k}{\partial x^\beta}\right)
\qe
for $f:M\times [0,\infty)\to N$ with initial values $f(x,0)=g(x)$. A solution is shown to exist for all times $0\le t < \infty$ under the assumption that $N$ has nonpositive sectional curvature and to converge to a solution of (\rf{13}) for $t\to \infty$ under the geometric assumptions of the theorem or the corollaries. In order to see the relevance of these assumptions, let us consider the following example:\\
On $\R^2$, consider the affine transformations
\bel
\label{ex1}
(x,y)\to (x+ny +m +\frac{1}{2} n^2, y+n)
\qe
for $m,n \in \Z$. The quotient of $R^2$ by this action of $\Z^2$ then is a compact affine manifold $M$, see e.g. \cite{GH}. 
\bel
\label{ex2}
{\tilde g}:\R^2 \to \R^1, (x,y)\mapsto x-\frac{1}{2}y^2
\qe
then is a map which equivariant w.r.t. the homomorphism $(m,n)\to m$ (i.e., ${\tilde g}(x+ny +m +\frac{1}{2} n^2, y+n)={\tilde g}(x,y)+m$  and therefore induces a map
\bel
\label{ex3}
g:M \to S^1
\qe
where $S^1=\R^1/\Z$. We consider the heat flow on $\R^2$,
\bel
\label{ex4} 
\frac{\partial \phi}{\partial t} =\Delta \phi
\qe
with initial values $\phi(x,y,0)={\tilde g}(x,y)$ where $\Delta$ is the standard Laplace operator $\displaystyle \frac{\partial^2}{\partial x^2}+\frac{\partial^2}{\partial y^2}$. The solution of (\rf{ex4}) is given by 
\bel
\label{ex5}
\phi(x,y,t)= x-\frac{1}{2}y^2-t
\qe
and therefore, it stays equivariant for all $t>0$. For $t\to \infty$, it disappears at infinity and does not converge to a harmonic function.\\
This is not precisely the situation considered here because the Laplace operator $\Delta$ is not invariant under the action  of $\Z^2$ on $\R^2$, but since the solution $\phi$ nevertheless stays equivariant, this does not matter. (Actually, an invariant metric is given by 
\bel
\label{ex6}
\gamma_{\alpha \beta}(x,y) = \left[
                 \begin{array}{cc}
                   1 & -y \\
                  -y & y^2 + 1 \\
                 \end{array}
               \right]
\qe
which is not K\"ahler affine.)

\section{Proof of the main theorem}\label{section4}
We shall abbreviate (\rf{15}) as
\bel
\label{20}
\frac{\partial f}{\partial t}=\sigma(f).
\qe
Since this is a system of parabolic differential equations, the existence of a solution on a
short time interval $[0,\tau)$ and, more generally, the openness of
the existence interval follow from general results about parabolic
equations. The first difficult step of the proof will now consist in
showing the closedness of the existence interval. For that step, we
shall need the nonpositive sectional curvature of the target. The
second  step will then be to show that the solution of (\rf{20})
converges to an affine harmonic map as $t\to \infty$. For that, we
need to show in particular that $f_t \to 0$ as $t\to \infty$. For that
step, we shall need to use the homotopic nontriviality condition in
addition to the nonpositive sectional curvature.\\
We now carry out the first step. It will be divided into several
substeps.
\begin{enumerate}
\item Let $f(x,t,s)$ be a family of solutions of (\rf{20}) depending on a parameter
  $s$. We then compute, using (\rf{20}) to convert third derivatives
  into curvature terms by the standard commutation formula for
  covariant derivatives
\bel
\label{21}
\left( \gamma^{\delta \epsilon} \frac{\partial^2}{\partial x^\delta
    \partial x^\epsilon}- \frac{\partial}{\partial t} \right) \left(
  g_{ij} \frac{\partial f^i}{\partial s} \frac{\partial f^j}{\partial
    s} \right)=2\gamma^{\delta \epsilon}\left( g_{ij} \frac{\partial^2
    f^i}{\partial x^\delta \partial s} \frac{\partial^2 f^j}{\partial
    x^\epsilon \partial
    s} - \frac{1}{2} R_{ijkl}\frac{\partial f^i}{\partial
    s}\frac{\partial f^j}{\partial x^\delta}\frac{\partial
    f^k}{\partial s}\frac{\partial f^l}{\partial x^\epsilon}\right)
\qe
where $R_{ijkl}$ is the curvature tensor of the target manifold
$N$. (A more detailed computation will be given in the next step.) Since we assume that the latter has nonpositive sectional
curvature, we conclude
\bel
\label{22}
\left( \gamma^{\delta \epsilon} \frac{\partial^2}{\partial x^\delta
    \partial x^\epsilon}- \frac{\partial}{\partial t} \right) \left(
  g_{ij} \frac{\partial f^i}{\partial s} \frac{\partial f^j}{\partial
    s} \right)\ge 0.
\qe
One such family of solutions is obtained by a time shift,
\bel
\label{23}
f(x,t,s):=f(x,t+s).
\qe
We use this to obtain
\begin{lemma}
\label{lem1}
\bel
\label{24}
\sup_{x\in M} g_{ij} \frac{\partial f^i}{\partial t} \frac{\partial
  f^j}{\partial t}
\qe
is nonincreasing in $t$ for a solution of (\rf{20}).
\end{lemma} 
\begin{proof}
Applying (\rf{22}) to the family (\rf{23}) yields
\bel
\label{25}
\left( \gamma^{\delta \epsilon} \frac{\partial^2}{\partial x^\delta
    \partial x^\epsilon}- \frac{\partial}{\partial t} \right) \left(
  g_{ij} \frac{\partial f^i}{\partial t} \frac{\partial f^j}{\partial
    t} \right)\ge 0,
\qe
and the maximum principle for subsolutions of parabolic equations then
implies the result.
\end{proof}

\item We consider
\bel
\label{31}
\eta(f):= \gamma^{\alpha \beta}g_{ij}(f(x,t))\frac{\partial
  f^i}{\partial x^\alpha}\frac{\partial
  f^j}{\partial x^\beta}
\qe
As in (\rf{22}), we want to apply the operator $\displaystyle \gamma^{\delta \epsilon} \frac{\partial^2}{\partial x^\delta
    \partial x^\epsilon} \frac{\partial}{\partial t}$ to this
  expression. This time, however, we also have to deal with
  derivatives of the domain metric. In order to simplify the
  computation, we use the standard device of orthonormal frames at the
  point under consideration. For the target, we may assume
  $g_{ij}=\delta_{ij}, g_{ij,k}=0$. For the domain, we may also assume
$\gamma_{\alpha \beta}=\delta_{\alpha \beta}$, but not necessarily
also the vanishing of the first derivatives. We then compute, using
subscripts for partial derivatives, 
\ba
\nonumber
& &\left(  \frac{\partial^2}{\partial x^\delta
    \partial x^\delta} \frac{\partial}{\partial t}\right) \eta(f)\\
\nonumber
&=& f^i_{x^\alpha x^\delta} f^i_{x^\alpha x^\delta}\\
\nonumber
&+& {\gamma^{\alpha \beta}}_{,\delta} (f^i_{x^\alpha x^\delta}f^i_{x^\beta}
+f^i_{x^\alpha }f^i_{x^\beta x^\delta}) +{\gamma^{\alpha
    \beta}}_{,\delta \delta} f^i_{x^\alpha}f^i_{x^\beta}\\
\label{32}
&-& R_{ijkl}f^i_{x^\alpha} f^j_{x^\delta}f^k_{x^\alpha} f^l_{x^\delta}
\ea
where we have again used the equation (\rf{20}). Using the Schwarz
inequality to handle the terms with first derivatives of the domain metric, the nonpositivity of the curvature of $N$ and rewriting
the result in general coordinates, we therefore obtain
\bel
\label{33}
\left( \gamma^{\delta \epsilon} \frac{\partial^2}{\partial x^\delta
    \partial x^\epsilon} \frac{\partial}{\partial t} \right) \eta(f)
\ge -c \eta(f) +\frac{1}{2} |D^2 f(.,t)|^2
\qe
with some constant $c$. In particular, $\eta(f)$ satisfies a linear
differential inequality, and we therefore obtain
\bel
\label{34}
\eta(f(x,t)) \le c \sup_{t_0 \le \tau \le t} \int_M \eta(f(.,\tau),
\qe
for any $t_0>0$, see e.g. \cite{J0}, Section 3.3.  

\item Next, as in \cite{JY}, using Jacobi field estimates (see e.g. \cite{J0}, Section
  2.5 and in particular formula (2.5.6) and the one preceding it), we
  obtain 
\bel
\label{35}
\int_M \eta(f(.,t))\le c\int_M (\tilde{d}^2(f(.,t),f^0)-\inf_{z\in
  M}\tilde{d}^2(f(z,t),f^0(z)))\quad  + c
\qe
where $\tilde{d}(f(.,t),f^0(.))$ is the homotopy distance between the
initial map $f^0=f(.,0)$ and the map $f(.,t)$ at time $t$; the
homotopy distance  $\tilde{d}(f(x,t),f^0(x))$ for these two homotopic
maps is given by the length of the shortest geodesic from $f(x,t)$ to
$f^0(x)$ in the homotopy class of curves determined by the homotopy
between the maps. \\
Also, these Jacobi field estimates yield 
\bel
\label{36}
\gamma^{\alpha \beta} \frac{\partial^2}{\partial x^\alpha \partial
  x^\beta}\tilde{d}^2(f(.,t),f^0)\ge -c \tilde{d}(f(.,t),f^0).
\qe

\item We can now complete the first step and prove long time existence
  of a solution of (\rf{20}). \\
Since in (\rf{34}), we take a supremum over different times $\tau$, we
have to control the behavior of our solution at different times
against each other. We have by the triangle inequality
\bel
\label{37}
\tilde{d}^2(f(.,\tau),f^0)\le
2\tilde{d}^2(f(.,t),f(.,\tau))+2\tilde{d}^2(f(.,t),f^0).
\qe
Also,
\bel
\label{38}
\tilde{d}^2(f(.,t),f(.,\tau))\le |t-\tau| \sup_{\tau \le \sigma \le
  t}|f_t(.,\sigma)| \le c|t-\tau|
\qe
where the last inequality follows from Lemma \ref{lem1}. With these
inequalities at hand, we can use (\rf{34}) and (\rf{35}) to obtain for the
norm of the first derivative $df$ w.r.t. the spatial variable $x$
\bel
\label{39}
|df(x,t)|\le c \left( \int_M \tilde{d}^2(f(.,\tau),f^0) \right)^{1/2}
+c
\qe 
and from this then also
\bel
\label{40}
|df(x,t)|\le c \sup_{y\in M}\tilde{d}(f(y,\tau),f^0(y)) 
+c. 
\qe  
Using (\rf{38}) then yields
\bel
\label{41}
|df(x,t)|\le c (1+t).
\qe
(\rf{41}) and Lemma \ref{lem1} yield $C^1$-bounds for our solution
$f(x,t)$ of (\rf{20}). We thus look at (\rf{15}) as an inhomogeneous
linear parabolic system with bounded right hand side. We can then apply the regularity theory for
solutions of linear parabolic equations to get $C^{2,\alpha}$-bounds
by the standard bootstrapping argument. Such bounds then imply
closedness of the interval of existence, hence global existence. Thus,
we have shown
\begin{lemma}
\label{lem2}
For a target manifold $N$ of nonpositive sectional curvature, the
solution $f(x,t)$ of (\rf{20}) exists for all $t\ge 0$. 
\end{lemma}

\end{enumerate}

We now turn to the second step of the proof, the convergence of the
solution $f(x,t)$ of (\rf{20}) to an affine harmonic map for $t\to
\infty$. Here, we need to use the assumption of topological
nontriviality  as expressed in our theorem in addition to nonpositive
target curvature (for necessary background material on nonpositive
curvature, we may refer to, e.g., \cite{J0a}). Again, we divide the reasoning into several
substeps. 
\begin{enumerate}
\item Let $x_0 \in M$ be a point where $\tilde{d}(f(y,\tau),f^0(y))$
  attains its minimum. Using (\rf{36}) and applying the maximum
  principle on both the ball $B(x_0,R)$ of radius $R$ about $x_0$ and
  on its complement $M{\backslash}B(x_0,R)$, we obtain
\bel
\label{50}
\sup_{y\in M}\tilde{d}^2(f(y,\tau),f^0(y))\le  \sup_{z\in \partial
  B(x_0,R)}\tilde{d}^2(f(z,\tau),f^0(z)) + c(R) \sup_{y\in
  M}\tilde{d}(f(y,\tau),f^0(y))
\qe
where the constant $c(R)$ depends on the radius $R$. The boundary term
can be controlled as follows
\ba
\label{51}
& &\sup_{z\in \partial
  B(x_0,R)}\tilde{d}^2(f(z,\tau),f^0(z))\\
\nonumber
&\le&
\tilde{d}^2(f(x_0,\tau),f^0(x_0)) +2R \sup_{y\in
  M}\tilde{d}(f(y,\tau),f^0(y))(|df(y,t)|+|df^0(y)|)
\ea
Using (\rf{40}), (\rf{50}) and (\rf{51}), we obtain for a suitable choice of
  $R>0$
\bel
\label{52}
\sup_{y\in M}\tilde{d}^2(f(y,\tau),f^0(y))\le \inf_{y\in
  M}\tilde{d}^2(f(y,\tau),f^0(y)) +c \sup_{y\in
  M}\tilde{d}(f(y,\tau),f^0(y)).
\qe 

\item Combining (\rf{35}) and (\rf{52}),
\bel
\label{53}
\int_M \eta(f(.,t))\le c \sup_{y\in
  M}\tilde{d}(f(y,\tau),f^0(y)) +c.
\qe 
Using then (\rf{34}) gives the pointwise estimate
\bel
\label{54}
|df(x,t)|\le c (\sup_{y\in
  M}\tilde{d}(f(y,\tau),f^0(y)))^{1/2} +c.
\qe 
Therefore, for any $x_1,x_2 \in M$, letting $\tilde{f}$ denote the
lift to universal covers, 
\bel
\label{55}
d(\tilde{f}(x_1,t),\tilde{f}(x_2,t))\le c (\sup_{y\in
  M}\tilde{d}(f(y,\tau),f^0(y)))^{1/2} +c.
\qe 

\item The essential point of the proof will be to exclude that for
  some sequence $t_n \to \infty$ and for some, and by (\rf{52}) then for
  all, $y\in M$,
\bel
\label{56}
\tilde{d}(f(y,t_n),f^0(y)) \to \infty. 
\qe
For $x \in M$, we let $\gamma^n_x$ be the geodesic from
$f^0(x)$ to $f(x,t_n)$ in the right homotopy class, i.e., the one
determined by the homotopy between the maps $f^0$ and
$f(.,t_n)$. Their length $T_n$ will then go to infinity, if
(\rf{56}). In fact, while the length depends on $~i$, by (\rf{52}), this is
inessential. \\
Since $N$ has nonpositive sectional curvature, the distance
\bel
\label{57}
d(\gamma^n_{x_1}(\tau), \gamma^n_{x_2}(\tau))
\qe
is a convex function of $\tau$. Since by (\rf{55}), this distance grows
at most like $(T_n)^{1/2}$, it must be bounded. Therefore, the
geodesic rays $\gamma_{x_i}$ that are the limits of $\gamma^n_{x_i}$ for $n\to
\infty$ (perhaps after a selection of a subsequence) satisfy
\bel
\label{58}
\kappa(x_1,x_2,\tau):=d(\gamma_{x_1}(\tau), \gamma_{x_2}(\tau))\le d(\gamma_{x_1}(0), \gamma_{x_2}(0))
\qe
for all positive $\tau$. \\
There are then two possibilities: Either $\kappa$ is decreasing in $t$
or 
constant. In fact, we may always assume the latter, by the
following observation. Since (\rf{58}) holds for any two points $x_1,
x_2$, we then also conclude that
\bel
\label{59}
\eta(f(x,t))
\qe
is a nonincreasing function of $t$ for every $x$, and it has to
decrease for some $x$ unless $\kappa$ is constant in $t$ for any two
points. When, however, we choose our initial values $f^0$ as a
harmonic map, i.e., one that minimizes $\int_M \eta(f(.)$, then
$\eta(f(x,t)$ can only be a nonincreasing function of $t$    for each
$x$ if it is constant. \\
Now, when $\kappa(x_1,x_2,\tau)$ is a constant function of $\tau$, it
generates a flat strip, since $N$ has nonpositive sectional
curvature. 
\item Also, if $f^0$ is energy minimizing, then for any $t\ge 0$,  the map
  $f^t(x):=\gamma_x(t)$ is also energy minimizing, by the same
  reasoning. We shall now use these energy minimizing maps to track
  our sequence $f(.,t)$ and to get time independent estimates. \\
We take a sequence $t_n \to \infty$ as above and write $f^n$ in place
of $f^{t_n}$. From the preceding constructions we obtain, in case
(\rf{56}9 holds,
\bel
\label{60}
\tilde{d}(f(.,t_n),f^n)\le c (\tilde{d}(f(.,t_n),f^0))^{1/2} +c.
\qe
We wish to get rid of the first term on the right hand side, i.e., we
want $f^n$ to track $f(.,t_n)$ uniformly. That will then give us some
control on the first derivatives of those maps w.r.t. $x$. \\
We can repeat the construction with $f^n$ in place of $f^0$. We have
two possibilities. Either after finitely many steps, we find some
energy minimizing map $\hat{f}^n$ with
\bel
\label{61}
\tilde{d}(f(.,t_n),\hat{f}^n)\le c,
\qe
or we generate a new flat direction from strips between geodesics rays
of constant distance as above in each step. In that case, however,
after finitely many steps, we have exhausted all possible directions,
and $N$ must be flat. In that case, it is elementary to track
$f(.,t_n)$ also in the desired manner, and in fact, we are then
dealing with linear parabolic equations which is much easier than the
nonlinear case. Thus, in either case, we may assume (\rf{61}).  \\
We may then apply the reasoning leading to (\rf{40}) with the variable
map  $\hat{f}^n$ in place of $f^0$ to obtain
\bel
\label{62}
|df(x,t)|\le c \sup_{y\in M}\tilde{d}(f(y,\tau),\hat{f}^n(y)) 
+c \le c.  
\qe  
 
\item With Lemma \ref{lem1} and (\rf{62}), we have uniform estimates for all
  first derivatives of $f(x,t)$, i.e., estimates that do not depend on
  $t$. Linear elliptic parabolic regularity theory then also yields
  higher order estimates, and we can then find a sequence $t_n \to
  \infty$ for which $f(.,t_n)$ converges smoothly to some smooth map $f_\infty$
  in the right homotopy class. It remains to show that $f_\infty$ is
  affine harmonic.

\item
We recall (\rf{21}) for the family $f(x,t,s):=f(x,t+s)$, that is, 
\ba
\label{72}
\nonumber
& &\left( \gamma^{\delta \epsilon} \frac{\partial^2}{\partial x^\delta
    \partial x^\epsilon}- \frac{\partial}{\partial t} \right) \left(
  g_{ij} \frac{\partial f^i}{\partial t} \frac{\partial f^j}{\partial
    t} \right)\\
&=&2\gamma^{\delta \epsilon}\left( g_{ij} \frac{\partial^2
    f^i}{\partial x^\delta \partial t} \frac{\partial^2 f^j}{\partial
    x^\epsilon \partial
    t} - \frac{1}{2} R_{ijkl}\frac{\partial f^i}{\partial
    t}\frac{\partial f^j}{\partial x^\delta}\frac{\partial
    f^k}{\partial t}\frac{\partial f^l}{\partial x^\epsilon}\right).
\ea
Since we know from Lemma \ref{lem1} that $g_{ij} \frac{\partial f^i}{\partial t} \frac{\partial f^j}{\partial
    t}$ stays bounded in $t$, and since both terms on the right
hand side of (\rf{72}) are nonnegative, they both have to converge to 0
for $t\to \infty$. The asymptotic vanishing of the first term means
that
\bel
\label{73}
\frac{\partial f(x,t)}{\partial t}
\qe
converges to a parallel section $v(x)$ along $f_\infty$ for $t\to
\infty$. This, however, is excluded in the assumptions of our
theorem. Therefore, 
\bel
\label{74}
\frac{\partial f(x,t)}{\partial t} \to 0 \text{ for } t \to \infty.
\qe
Thus, in the limit $t\to \infty$, the temporal derivative disappears
in (\rf{20}), and the elliptic system that we want to solve remains. 
This, together with the smooth convergence of $f(.,t_n)$ to
$f_\infty$, 
shows that $f_\infty$ solves the elliptic system, i.e., it is affine
harmonic. This completes the proof of our main theorem. In fact, it is
not hard to show now that the solution $f(.,t)$ of the parabolic
system converges to the solution $f_\infty$ of the elliptic system as
$t\to \infty$. 

\end{enumerate}

{\bf Remarks:}
\begin{enumerate}
\item Naturally, one can also treat the Dirichlet problem for affine
  harmonic maps. Here, one could either use the method of \cite{JY} or
  the general approach developed by von Wahl \cite{vW1,vW2} for
  parabolic systems that does not need a variational
  structure. When Dirichlet boundary values are given, they prevent a
  solution from eternally moving around the target manifold. Thus, the
  main problem that we had to overcome in the proof of our main
  theorem and for which we needed an additional topological assumption
  besides the geometric condition of nonpositive curvature is not
  present in the Dirichlet boundary value problem. Of course, boundary
  regularity then is an issue that needs treatment, but this can be
  achieved by the methods of the aforementioned papers. 
\item It should be possible and of interest in affine geometry to extend the method of Grunau and
  K\"uhnel \cite{GK} to show the existence of affine harmonic maps
  from a complete affine to a complete Riemannian manifold.
\end{enumerate}

\bigskip
{\bf Acknowledgements}
The second author is grateful to J. Jost for posing the problem and for
stimulating discussions.  The second author  was  supported by the The Scientific and Technologic Research Council 
of Turkey, 2219 Fellowship and Max Planck Institute of Mathematics in the Sciences.

\bibliographystyle{plain}

\bigskip
\noindent
\parbox[t]{.48\textwidth}{
J\"urgen  Jost\\
Max Planck Institute for Mathematics in the Sciences\\
Inselstrasse 22\\
D-04103 Leipzig, Germany\\
jost@mis.mpg.de} \hfill
\parbox[t]{.48\textwidth}{
Fatma Muazzez \c{S}im\c{s}ir\\
TOBB University of Economics and Technology\\
Department of Mathematics\\
S\"o\u{g}\"ut\"oz\"u Caddesi No: 43\\
TR-06560 Ankara, Turkey\\
msimsir@etu.edu.tr}
\end{document}